\documentclass[a4paper]{article}
\usepackage[english]{babel}

\title{Decomposition of fuzzy relations: \\An application to fuzzy preferences}
\author{María Jesús Campión, Esteban Induráin \& Armajac Ravent\'os-Pujol}
\date{}

\usepackage{hyperref}

\usepackage{amsfonts,amssymb,amsthm}
\usepackage{amsmath}
\usepackage[T1]{fontenc}
\usepackage{graphicx}
\usepackage{float}
\usepackage{diagbox}
\usepackage{makecell}
\usepackage{lscape}
\usepackage{subcaption}
\usepackage{comment}
\usepackage{chngcntr}
\usepackage{apptools}

\usepackage{pgf,tikz, tikz-3dplot}
\usepackage{pgfplots}

\usetikzlibrary{arrows}
\usetikzlibrary{patterns,intersections}
\usetikzlibrary{math}
\usetikzlibrary{decorations.pathreplacing,angles,quotes}
\usetikzlibrary{decorations.markings}
\usetikzlibrary{positioning}
\usetikzlibrary{through}
\usetikzlibrary{tikzmark}
\usetikzlibrary{angles,quotes}
\usetikzlibrary{calc}
\pgfplotsset{compat=1.17}
\usepgfplotslibrary{fillbetween}

\numberwithin{equation}{section}
\theoremstyle{plain}
\newtheorem{theorem}{Theorem}[section]
\newtheorem{proposition}[theorem]{Proposition}

\theoremstyle{definition}
\newtheorem{definition}[theorem]{Definition}
\newtheorem{example}[theorem]{Example}

\theoremstyle{remark}
\newtheorem{remark}[theorem]{Remark}

\begin{document}
\maketitle

\begin{abstract}
In this article, working in the spirit of the classical Arrovian models in the fuzzy setting and their possible extensions, we go deeper into the study of some type of decompositions defined by t-norms and t-conorms. This allows us to achieve characterizations of existence and uniqueness for such decompositions and we provide rules to obtain them under some specific conditions. We conclude by applying such results to the study of fuzzy preferences, a key concept in Fuzzy Social Choice.
\end{abstract}

\section{Introduction}
Let us consider an individual that defines her/his preferences on a nonempty set $X$. To start with, let us assume that those preferences are crisp binary relations on $X$, so that given two elements $x,y \in X$ and a binary relation $\mathcal{R}$ either $x$ is related to $y$ through $\mathcal{R}$ or it is not. In other words, understanding $\mathcal{R}$ as a subset of the Cartesian product $X \times X$ we have that either $(x,y) \in \mathcal{R}$ or else $(x,y) \notin \mathcal{R}$. In that crisp setting, $\mathcal{R}$ is then a crisp subset of $X \times X$, that is the characteristic function $\chi_{\mathcal{R}}$ that defines $\mathcal{R}$ can only take the values $0$ or $1$. Thus $\chi_{\mathcal{R}}(x,y)= 1$ means that $x$ is related to $y$ by means of $\mathcal{R}$, whereas $\chi_{\mathcal{R}}(x,y)= 0$ means that $x$ is not related to $y$ through $\mathcal{R}$.

In the crisp approach a \emph{preference} on $X$ is usually understood as a binary relation that is transitive and complete\footnote{A complete binary relation is also known in the literature as a total binary relation. We also remark that a complete binary relation is also reflexive.}. This is also known as a \emph{total preorder}. We will use the notation $\succsim$. Associated to it we may consider also two binary relations, namely $\succ$ and $\sim$, respectively called the \emph{strict preference} and the \emph{indifference}, and defined by $x \succ y$ if and only if it is not true that $y \succsim x$, and $x \sim y$ if and only if both $x \succsim y$ and $y \succsim x$ are true ($x,y \in X)$. With this definition, the total preorder $\succsim$ decomposes as $\succ \cup \sim$, and this decomposition is unique (see Section \ref{sec:preliminaries}). 

Unlike the crisp approach, in the fuzzy setting a \emph{fuzzy binary relation} $R$ on a nonempty set $X$ is now understood as a fuzzy set of the Cartesian product $X \times X$, so that now the \emph{indicator} (also know as the \emph{membership function}) of $R$ is a map $\chi_R: X\times X \rightarrow [0,1]$. It may  take any possible value from $0$ (no relationship at all, absolute negation) to $1$ (total relationship, total evidence). Intermediate values give us, then, an idea of uncertainty. 

To start now working in the fuzzy setting, some rather important facts should be pointed out:
\begin{itemize} 
    \item[1.] There are many non-equivalent extensions of the notions of completeness and/or transitivity, so that many possible definitions of the concept of a \emph{fuzzy total preorder} could appear (see Section \ref{sec:preliminaries}).
    \item[2.] Once a particular definition of a fuzzy total preorder, say $R$,  has been chosen, we may try to decompose it by means of two fuzzy binary relations $P$ and $I$ so that, in a way, the triplet $(R,P,I)$ behaves as the decomposition $(\succsim,\succ,\sim)$ of a crisp total preorder $\succsim$. However, the definitions of what we mean for a fuzzy decomposition are neither unique nor equivalent (see Section \ref{sec:generalDecompositions}).
    \item[3.] Again, even if we have chosen some particular definition of a fuzzy decomposition, it may happen that a given fuzzy preorder could still have more than one possible decompositions of that kind. In other words, now the decompositions are not necessarily unique (see Sections \ref{sec:strongDec} and \ref{sec:weakDec}).
\end{itemize}

The problem of choosing definitions when generalizing crisp concepts is widespread in the fuzzy literature \cite{billot2012economic, DeBaets1995333, Gibilisco2014, MONTERO198645, Richardson1998359}. Particularly, the study of preference structures is not an exception. For example, whereas some authors study decompositions into two fuzzy binary relations $P$ and $I$ \cite{Fono2005372,Llamazares2005, Richardson1998359}, others do it into three relations: $P$, $I$ and an additional one representing the idea of incomparability \cite{Diaz20082221, Diaz20101, FODOR199247, LLAMAZARES2003217, OVCHINNIKOV199215}.

Our study on decompositions started motivated by the need for a solid framework for fuzzy Arrovian models (see a complete survey in \cite[Chapter 4]{Gibilisco2014}). These models are a fuzzy extension of the well-known model introduced by Arrow in 1951, in the crisp setting, studying ranking aggregation \cite{MR0039976}. In these fuzzy models, the decomposition of a preference plays an essential role, and it can be decisive in determining when it is possible to aggregate preferences\footnote{We can illustrate the importance of decompositions comparing the models studied in \cite[Proposition 3.14]{Dutta1987215} and \cite[Proposition 3.1]{Banerjee1994121}. Their definition only differs in the chosen decomposition rule. In the former case, it is possible to aggregate preferences, but it is not possible in the second case.}. In most of this literature, decompositions are defined using t-norms and t-conorms (as generalizations of the crisp intersection and crisp union). Whereas in the earlier models, decompositions were defined using a t-norm and a t-conorm playing the role of the intersection and the union \cite{Dutta1987215}, in subsequent works, the requirements for decompositions were weakened: the authors only used a single t-conorm in addition to a condition to assure that the fuzzy decomposition was a generalization of the dichotomic one with values in $\{0,1\}$ \cite{Fodor, Fono2005372, Gibilisco2014, Richardson1998359}. One of our goals in our study is understanding the differences between both types of decomposition (in this article strong and weak decompositions) and the reason behind this change. 

The main difference in preference modelling in these Arrovian models and some of the works cited above as \cite{Diaz20082221, Diaz20101, FODOR199247, LLAMAZARES2003217, Llamazares2005, OVCHINNIKOV199215} is the definition of asymmetry for the strict preference $P$ as well as the nonexistence of an incomparability relation. The second one seems reasonable since the preferences are complete in the original Arrow's model. That is, two alternatives are always comparable. However, the choice of the asymmetry on fuzzy relations requires a deeper justification. Despite the most intuitive generalization requires $T(P(x,y), P(y,x))=0$ for some t-norm $T$ (see \cite{Llamazares2005}), most of the researchers in this area have stood for strict preferences without both conjugated pairs with positive degree (i.e. $P(x,y), P(y,x)>0$) and, as Dutta explained in one of the earliest articles \cite{Dutta1987215}, it is the most restrictive conception of asymmetry.

The reasons behind both choices make sense when we take into account the historical perspective. The introduction of fuzzy preferences in the Arrovian model was motivated by escaping from the Arrow's impossibility theorem \cite{MR0039976}. As we have explained before, when fuzziness is introduced, many possible extensions of the concepts and definitions (unique in the original model) arise. These new definitions must make sense in the problems social scientists are studying. We look for models allowing preference aggregation rules, but in many cases, it is needed to remain as much as possible close to a dichotomic interpretation. For instance, Duddy \& Piggins studied aggregation rules where individual preferences were rankings whereas social preferences were fuzzy \cite{Duddy20181}. Billot also argues that any transitivity defined by t-norms will contradict the \emph{independence of irrelevant alternatives principle}\footnote{The independence of irrelevant alternatives is a (more philosophic than mathematical) assumption that states that if you take as input a set of alternatives and their characteristics, the input should not depend on extra-group characteristics.} and they should be avoided \cite[Section 1.1.2]{billot2012economic}. In the same spirit, Dutta recognizes he chose the most restrictive form of asymmetry \cite{Dutta1987215} because his study was looking for the minimum fuzzification needed to obtain aggregation rules. Moreover, his assumption has not been the object of remarkable criticism in subsequent studies (e.g. \cite{Banerjee1994121, Fono2005372, Gibilisco2014, Richardson1998359}).\\
For example, consider the case in which $R(x,y)=0.9$ and $R(y,x)=0.05$. Here $x$ is clearly preferred over $y$, so under the assumption made by Dutta and others, it makes sense to state that $P(y,x)=0$. Moreover, compared with the situation $R'(x,y)=1$ and $R'(y,x)=0$, it makes sense that $P'(x,y)>P(x,y)$ because of the differences of the intensities between $R$ and $R'$. Then, we could argue that $P$ must be a fuzzy relation satisfying the most restrictive asymmetry when an alternative is clearly preferred over the other. However, this argument does not lack criticism: consider the case in which $R(x,y)=0.65$ and $R(y,x)=0.64$. Can we state that $x$ is clearly preferred over $y$? In this case, saying that $P(x,y)>0$ and $P(y,x)=0$ lacks real meaning, and it should not be taken into account by practitioners. Nevertheless, cases such as the first example have prevailed over the second, and the restrictive notion of asymmetry is imposed. One of its reasons is, probably, that \emph{Pareto condition}\footnote{Pareto condition states that if all agents strictly prefer an alternative over another, the social preference must do it (see e.g. \cite{Gibilisco2014}). The Pareto condition is the most accepted condition in the Arrovian model and is the backbone of the aggregation functions structure.} lose its importance if other types of asymmetries are accepted.

It is worth mentioning that not all social scientists share the same opinion. Llamazares, in his study of decompositions \cite{Llamazares2005}, justifies his choice of asymmetry quoting Blin \cite{Blin1974} and Barret \& Pattanai \cite{Barrett1989} \guillemotleft[the vagueness] \emph{arises through the multiplicity of dimensions underlying preferences}\guillemotright , and for that reason $P(x,y)>0$ and $P(y,x)>0$ can coexist. However, we have to mention that most of the works of Barret and Pattanaik in fuzzy preferences aggregation only work with strict preferences, and they do not have to deal with the decomposition problem nor the differences between weak and strict preferences \cite{Barrett19861, BARRETT19929}.

Having in mind the previous discussion, the scheme of the manuscript goes as follows:
 
After the introduction, a section of preliminaries furnishes the basic definitions that will be used along the whole paper. We pay attention there to classical concepts in fuzzy set theory.  Section \ref{sec:generalDecompositions} introduces then the key concept of a decomposition of a fuzzy binary relation. Among the main families of decompositions, we will consider the so-called strong and weak ones, respectively analyzed in Sections \ref{sec:strongDec} and \ref{sec:weakDec}. In Section \ref{sec:decompositionPreferences}, we apply the results of the previous sections on the decomposition of preferences in the context of fuzzy Arrovian models. Finally, in \ref{sec:restrictDomains}, we will give some geometrical ideas about the existence of decompositions when spaces of fuzzy binary relations are restricted by transitivity or completeness conditions. And we close the paper in \ref{app:examples} applying our results to some of the most known t-norms and t-conotms.

\section{Preliminaries}
\label{sec:preliminaries}

In this article, $X$ will always denote a nonempty set.

\begin{definition} A binary relation $\mathcal{R}$ on $X$ is a subset of the Cartesian product $X \times X$. Given two elements $x,y \in X$, we will use the standard notation $x \mathcal{R} y$ to express that the pair $(x,y)$ belongs to $\mathcal{R}$.

Associated to a binary relation $\mathcal{R}$ on $X$, we consider its \emph{negation} \rm (respectively, its \emph{transpose}\rm) as the binary relation $\mathcal{R}^c$ (respectively, $\mathcal{R}^t$) on $X$ given by $(x,y) \in \mathcal{R}^c  \Leftrightarrow (x,y) \notin \mathcal{R}$ for every $ x,y \in X$ (respectively, given by $(x,y) \in \mathcal{R}^t  \Leftrightarrow (y,x) \in \mathcal{R}, \ $ for every $ x,y \in X)$. We also define the \emph{adjoint} \rm $\mathcal{R}^a$ of the given relation $\mathcal{R}$, as $\mathcal{R}^a = (\mathcal{R}^t)^c$. 

A binary relation $\mathcal{R}$ defined on a set $X$ is called: 
\begin{itemize}
   \item[(i)] \emph{reflexive} \rm if $x \mathcal{R} x$ holds for every $x \in X$,
   \item[(ii)] \emph{irreflexive} \rm if $\neg(x \mathcal{R} x)$ holds for every $x \in X$, 
   \item[(iii)] \emph{symmetric} \rm if $\mathcal{R}$ and $\mathcal{R}^t$ coincide,
   \item[(iv)] \emph{antisymmetric} \rm if $\mathcal{R} \cap \mathcal{R}^t \subseteq   \{ (x,x): x \in X \}$,
   \item[(v)] \emph{asymmetric} \rm if $\mathcal{R} \cap \mathcal{R}^t = \varnothing$,
   \item[(vi)] \emph{complete} \rm if $\mathcal{R} \cup \mathcal{R}^t = X \times X$,
   \item[(vii)] \emph{transitive} \rm if $x \mathcal{R} y \wedge y \mathcal{R} z \Rightarrow x \mathcal{R} z $ for every $x,y,z \in X$.
 \end{itemize}
 \end{definition} 

In the particular case of a set $X$ where some kind of \emph{ordering} \rm has been defined, the standard notation is different. 
\begin{definition} A \emph{preorder} $\succsim$ on $X$ is a binary relation on $X$ which is reflexive and transitive.
An antisymmetric preorder is said to be an \emph{order}. A \emph{total preorder} \rm $\succsim$ on a set $X$ is a preorder such that if $x,y \in X$ then $(x \succsim y) \vee (y \succsim x)$ holds.  
If $\succsim$ is a preorder on $X$, then as usual we denote the associated \emph{asymmetric} relation by $\succ$ and the
associated \emph{equivalence} relation by $\sim$
and these are defined by $x \succ y  \Leftrightarrow (x \succsim y) \wedge\neg(y \succsim x)$ and $x \sim y \Leftrightarrow (x\succsim y) \wedge (y \succsim x)$. 
 
A total preorder $\succsim$ defined on $X$ is usually called a (crisp) \emph{preference} on $X$. 
\end{definition}

\begin{definition} Let $\mathcal{R}$ be a binary relation on $X$. We say that $\mathcal{R}$ decomposes into a symmetric binary relation (denoted by $\mathcal{R}_s$) and an asymmetric binary relation (denoted by $\mathcal{R}_a$) if $\mathcal{R}=\mathcal{R}_a\cup \mathcal{R}_s$ and $\mathcal{R}_a\cap \mathcal{R}_s=\emptyset$.
\label{def:crispDec}
\end{definition}

\begin{proposition} Let $\succsim$ be a preorder on a set $X$. It decomposes into $\succ$ and $\sim$. Besides,  it is the unique decomposition of $\succsim$ into an asymmetric and a symmetric relations. 
\begin{proof}
It is straightforward to see that $\succsim$ decomposes into $\succ$ and $\sim$, and as a matter of fact  $\succ$ is asymmetric and $\sim$ symmetric. It only remains to see that it is the unique one. Suppose that $\succsim$ decomposes into two relations $\mathcal{R}_a$ and $\mathcal{R}_s$ with $\mathcal{R}_a$ asymmetric and $\mathcal{R}_s$ symmetric. First we will prove that $\succ = \mathcal{R}_a$. \\
Suppose that $(x,y)\in \succ$ but $(x,y)\notin \mathcal{R}_a$, then as $(x,y)\in \succsim =\mathcal{R}_a\cup \mathcal{R}_s$ we can assure that $(x,y)\in\mathcal{R}_s$. Now, as $\mathcal{R}_s$ is symmetric, then it holds that $(y,x)\in\mathcal{R}_s$, which implies that $(y,x)\in\succsim$, arriving at $(x,y)\in\sim$, a contradiction. We conclude that $\succ\subseteq \mathcal{R}_a$.\\
If we suppose that $(x,y)\in \mathcal{R}_a$ but $(x,y)\notin\succ$, by the same argument as before we obtain that $(x,y)\in\sim$, because $(x,y)\in \succsim =\succ\cup \sim$. By asymmetry of $\mathcal{R}_a$ we can assure that $(y,x)\in \mathcal{R}_a$ does not hold, but $(y,x)\in \sim\subseteq \succsim$, so $(y,x)\in \mathcal{R}_s$. Finally, as $\mathcal{R}_s$ is symmetric, $(x,y)\in \mathcal{R}_a\cap \mathcal{R}_s$, a contradiction because this intersection has to be empty. So $\mathcal{R}_a\subseteq \succ$, and $\succ=\mathcal{R}_a$.\\
If $\mathcal{R}_a=\succ$ is straightforward to see that $\mathcal{R}_s=\sim$.
\end{proof}
\end{proposition}

\begin{definition} A \emph{fuzzy subset} $H$ of $X$ is defined as a function $\mu_H: X \rightarrow [0,1]$.  The function $\mu_H$ is called the \emph{membership function} \rm  of $H$.  In the particular case when $\mu_H$ is dichotomic and takes values in $\lbrace 0, 1 \rbrace$, the corresponding subset defined by means of  $\mu_H$ is a subset of $X$ in the classical crisp sense\footnote{The term \emph{crisp} is usually understood in these contexts as meaning non-fuzzy.}.
\end{definition}

\begin{definition} A \emph{fuzzy binary relation on $X$} is a function $R: X\times X \rightarrow [0,1]$. We say that $R$ is \emph{symmetric} if $R(x,y)=R(y,x)$ for every $x,y\in X$ and we say that $R$ is \emph{asymmetric} if for every $x,y\in X$, $R(x,y)>0$ implies that $R(y,x)=0$ holds. Moreover, $R$ is said to be \emph{reflexive} if $R(x,x)=1$ for any $x\in X$.
\end{definition} 

\begin{definition} \label{tria} A \emph{triangular norm} (\emph{t-norm} for short) is a function $\allowbreak T:[0,1]\times[0,1]\longrightarrow [0,1]$ satisfying the following properties: 
	\begin{itemize}
		\item [i)] Boundary conditions: $T(x,0) = T(0,x) = 0$, and $T(x,1) = T(1,x) = x$,  for every $x \in [0,1]$.  
		\item [ii)] Monotonicity: $T$ is non-decreasing with respect to each variable, that is $ $ if $x_1\le x_2$ and $y_1\le y_2$, then  $T(x_1, y_1) \le T(x_2, y_2)$ holds true.
		\item [iii)] $T$ is commutative:  $T(x,y) = T(y,x)$ holds for every $x,y\in [0,1]$.
		\item [iv)] $T$ is associative:  $T(x , T (y,z))= T(T(x,y),z)$ holds for any $x,y,z\in [0,1]$.
	\end{itemize} 
\end{definition}

\begin{definition} A \emph{triangular conorm} (\emph{t-conorm} for short) is a function $S :[0,1]\times[0,1]\longrightarrow [0,1]$ satisfying the following properties:
	\begin{itemize}
		\item [i)] Boundary conditions: $S(x,0) = S(0,x) = x$, and $S(x,1) = S(1,x) = 1$,  for every $x \in [0,1]$. 
		\item [ii)] Monotonicity: $S$ is non-decreasing with respect to each variable, that is if $x_1\le x_2$ and $y_1\le y_2$, then  $S(x_1, y_1) \le S(x_2, y_2)$ holds true.
		\item [iii)] $S$ is commutative:  $S(x,y) = S(y,x)$ holds for every $x,y\in [0,1]$.
		\item [iv)] $S$ is associative:  $S(x , S (y,z))= S(S(x,y),z)$ holds for any $x,y,z\in [0,1]$.
	\end{itemize} 
	\label{tria2}
\end{definition}

Henceforth the symbol $T$ will denote a t-norm whereas the symbol $S$ will stand for  a t-conorm.

\begin{remark} In this article, we will use t-norms and t-conorms to define the decomposition of a fuzzy binary relation. However, in the proofs below, we will only use properties $i$ and $ii$ of Definitions \ref{tria} and \ref{tria2}. This fact shows that we could define a decomposition with respect to a more general category of operators (for example pre-aggregation functions \cite{7737700}, copulas \cite{Nelsen2006}, fusion functions \cite{9288714}, mixture functions \cite{MARQUESPEREIRA200343}, overlap functions \cite{DEMIGUEL201981}, penalty functions \cite{doi:10.1080/03081079708945181}, etc.). 
We have decided to use t-norms and t-conorms because they are the fuzzy generalizations of the classic intersection and union operators and, more importantly, because we want to be coherent with the preexisting literature that analyzes, studies and uses decompositions in practice \cite{Fono2005372, Gibilisco2014, Richardson1998359}.
\label{remark:MoreOperators}
\end{remark}

Let $R_1$, $R_2$ and $R_3$ be fuzzy binary relations on $X$. Then, the notation  $R_1=T(R_2,R_3)$ will mean from now on that, for every $x,y\in X$, $R_1(x,y)=T(R_2(x,y),R_3(x,y))$ holds. Obviously, we use a similar convention for a t-conorm $S$.

\begin{definition} Let $T$ and $S$ respectively be a t-norm and a t-conorm. An element $a\in (0,1)$ is said to be \emph{a $0$-divisor of $T$} if there exists some $b\in(0,1)$ such that $T(a,b)=0$. Also, an element $a\in (0,1)$ is called \emph{a $1$-divisor of $S$} if there exists some $b\in (0,1)$ such that $S(a,b)=1$.
\end{definition}

\begin{definition} A t-norm $T$ (respectively a t-conorm $S$) is called \emph{strict} if it is strictly increasing on each variable on $(0,1]^2$ (respectively $[0,1)^2$). In other words, for all $x,y,z\in[0,1]$ with $x<y$, $T(z,x)<T(z,y)$ if $0<z$ (respectively, $S(z,x)<S(z,y)$ if $z<1$). 
	\label{DefStrict}
\end{definition}

\section{Decompositions of fuzzy binary relations}
\label{sec:generalDecompositions}

In the present manuscript, we study decompositions based on t-norms and t-conorms respectively playing the crisp intersection and union role. This approach motivates the following definition: 

\begin{definition} Let $R,P$ and $I$ be fuzzy binary relations on $X$ (with $P$ an asymmetric fuzzy binary relation and $I$ a symmetric binary relation). Let $T$ (respectively, $S$) stand for a triangular norm (respectively, conorm).  We say that $R$ \emph{strongly decomposes},  with respect to $(T,S)$, into $(P,I)$, if $R=S(P,I)$ and $T(P,I)=0$.
\label{strongDec}
\end{definition}

Even if the decomposition defined above is the natural generalization of the usual crisp decomposition, it is not widespread in the literature. Indeed, many authors use some weaker versions \cite{Banerjee1994121, Dutta1987215, Fono2005372, Gibilisco2014, Richardson1998359}. The following definition follows the same spirit, but disregarding the role that any t-norm could play:

\begin{definition} Let $R$ be a fuzzy binary relation on a set $X$. Let $S$ be a t-conorm. We say that $R$ \emph{weakly decomposes} into $(P,I)$ (with $P$ an asymmetric fuzzy binary relation and $I$ a symmetric binary relation) if $R=S(P,I)$ and if for every $x,y\in X$ it holds true that if $I(x,y)=1$, then $P(x,y)=0$.
\label{weakDec}
\end{definition}

We have two comments about the definitions above. First, notice that the weak decomposition is more general than the stronger one. Moreover, they are not equivalent. The following example proves it:
\begin{example} Let $X=\{x,y\}$ and let $R$ be a fuzzy binary relation defined on $X$ as follows: $R(x,y)=R(x,x)=R(y,y)=1$ and $R(y,x)=0.5$. We define the symmetric relation $I$ as $I(x,y)=0.5$ and $I(x,x)=I(y,y)=1$. If we define $P$ as $P(x,y)=1$ and $0$ on the remaining cases, we can check that $R$ weakly decomposes into $(P,I)$ with respect the t-conorm maximum $S_{\max}$. However, $R$ does not strongly decomposes into $(P,I)$ with respect to the conorm $S_{\max}$ and any t-norm $T$ because $T(1,0.5)=0.5\neq 0$.  
\end{example}

The second observation is about the second condition in Definition \ref{weakDec}. It has been introduced after removing the t-norm from Definition \ref{strongDec} because we need decompositions being a generalization of the crisp one (Definition \ref{def:crispDec}). It implies that crisp relations must decompose into two crisp relations $(P,I)$ and these must be unique. But, for example, without this condition, if $X=\{x,y\}$ and a fuzzy binary relation $R$ on $X$ is defined as $R(x,y)=R(y,x)=R(x,x)=R(y,y)=1$, it decomposes to $(P,I)$ with $I(x,y)=I(x,x)=I(y,y)=1$ and $P(x,y)=0.5$ and $0$ otherwise.\\

Finally, in a previous paper \cite{Raventos-Pujol202003} the authors proved some necessary conditions for these kinds of decompositions. We will include them here for the sake of completeness.

\begin{proposition} Let $S$ be a t-conorm. Let  $R$ be a fuzzy binary relation defined on $X$. Assume that $R$ can be expressed as $R=S(P,I)$ (where $I$ (respectively, $P$) is a symmetric (respectively, asymmetric) fuzzy binary relation). Then, for every $x,y \in X$ it holds that $I(x,y)=\min\{R(x,y),R(y,x)\}$.
	\begin{proof} See the proof in \cite[Proposition 2]{Raventos-Pujol202003}.
	\end{proof}
	\label{simmin}
\end{proposition}

\begin{proposition} Let $S$ be a t-conorm. Every fuzzy relation $R$ on $X$ can be expressed as $R=S(P,I)$ (where $I$ (respectively, $P$) is a symmetric (respectively, asymmetric) fuzzy binary relation) if and only if $S$ is continuous in the first coordinate (with respect to the Euclidean topology on the unit interval $[0,1]$). 
	\begin{proof} See the proof in \cite[Proposition 3]{Raventos-Pujol202003}
	\end{proof}
	\label{Cdecomp}
\end{proposition}

The proof of Proposition \ref{Cdecomp} in \cite{Raventos-Pujol202003} gave use a method to obtain a candidate for the asymmetric component $P$ computed from $R$. In other words, if $R$ decomposes, one of its decompositions is the pair $(P,I)$ provided by this method.\\
As Proposition \ref{simmin} states, the $I$ component is a minimum, whereas the candidate for the $P$ component is defined as ${P(x,y)=\inf\{t\in [0,1]: S(t,I(x,y))\ge R(x,y)\}}$. To abbreviate that, we use $P(x,y)=I(x,y) \searrow_S R(x,y)$.

\begin{remark} Notice that since $S$ is commutative, continuity on one of its coordinates implies continuity on both coordinates (considered separately). However, it does not imply the continuity for a general symmetric real function. However, since $S$ is non-decreasing, both conditions are equivalent (see e.g. \cite[Proposition 1]{kruse69}). Then, in Proposition \ref{Cdecomp} being $S$ continuous in the second coordinate is equivalent to $S$ being continuous.\\
In spite of this fact, in the article we will specify over which component $S$ must be continuous. We proceed this way because the same results could be useful for the operators we have listed on Remark \ref{remark:MoreOperators}.
\end{remark}

In the literature, there are similar results to that stated in Proposition \ref{simmin} \cite{Fono2005372, Gibilisco2014, Richardson1998359}. In fact, even Fono and Andjiga \cite{Fono2005372} use a similar strategy than the one we use in Proposition \ref{Cdecomp}. However, the definitions of decompositions used in these works are more general than the weak decomposition defined in the present manuscript. Their definition is so general that they have to impose additional constraints. They end up studying \emph{regular decompositions}, and these are actually less general than Definition \ref{weakDec}.\\
In addition, in contrast with the work of Fono and Andjiga \cite{Fono2005372}, here, we do not assume a priori any continuity condition on $S$. Whereas they assume upper continuity in the definition of its decomposition, we have seen using Proposition \ref{Cdecomp} above that the continuity on the first component is a necessary condition for the existence of our decompositions. So, the continuity condition over $S$ is a consequence of $S$ admitting decompositions, and it is not imposed a priori.\\

In Section \ref{sec:strongDec} we will characterize the existence and the uniqueness of strong decompositions, whereas on Section \ref{sec:weakDec} we do the same for weak decompositions.\\
Given a fuzzy binary relation $R$, we could study the conditions that should satisfy $R$ in order to decompose. However, in the Social Choice models \cite{Banerjee1994121, Dutta1987215, Gibilisco2014, Richardson1998359} we usually work with it is needed that all the binary relations of a suitable set admit a decomposition. For that reason, we will start studying under which conditions all elements in $\mathcal{BR}_X$ (the set of fuzzy binary relations on $X$) decompose. Next, in \ref{sec:restrictDomains} we will explore the decomposition on subsets of $\mathcal{BR}_X$, mainly the subsets obtained by imposing some kind of transitivity and connectedness to fuzzy binary relations.

\section{Analysis and structure of strong decompositions of fuzzy binary relations}
\label{sec:strongDec}

This section contains two main results. We characterize the existence and the uniqueness of strong decompositions based upon  the divisors of a triangular norm $T$ and a t-conorm $S$. First, we start with some more specific (and preparatory) results which will motivate the general theorems to be proved next.\\

First, we focus on the case in which the used t-norm is the drastic t-norm $T_D$ (see its definition in Example \ref{example:standarnorms}). The following results prove that being decomposable with respect $T_D$ is a necessary condition for being decomposable with respect to any other t-norm.

\begin{proposition} Let $R$ be a fuzzy binary relation on $X$. Let $S$ be a t-conorm and $T$ a t-norm. If $R$ strongly decomposes into $(P,I)$ (with respect to $S$ and $T$), then $R$ decomposes into $(P,I)$ with respect to $S$ and $T_D$.
	\begin{proof} It is straightforward if we use the following fact: every t-norm is bounded by below by the drastic t-norm.
	\end{proof}
\end{proposition}

Considering the result above, we characterize as a first step the existence of decompositions with respect $T_D$.

\begin{proposition} Let $S$ be a t-conorm continuous in the first component. Any $R\in \mathcal{BR}_X$ is strongly decomposable with respect to $S$ and $T_D$ if and only if every $t\in (0,1)$ is a $1$-divisor (with respect to $S$).
	\begin{proof} 
	To prove the direct implication, let $t\in (0,1)$ and consider a fuzzy binary relation $R$ such that for a pair of values $x,y\in X$, $R(x,y)=1$ and $R(y,x)=t$. By hypothesis there exist $P$ and $I$ such that $1=R(x,y)=S(P(x,y),I(x,y))$ and $P(x,y)\neq 1$ because $T_D(P(x,y),I(x,y))=0$. Hence, $t=I(x,y)$ is a $1$-divisor.
	
	To prove the converse implication, it is enough to see that the weak decomposition $(P,I)$ proposed in the proof of Proposition \ref{Cdecomp} satisfies, for every $x,y\in X$, that $T_D(\allowbreak P(x,y),I(x,y))=0$. If we assume  that $T_D(P(x,y)\allowbreak ,I(x,y))\neq 0$, then $P(x,y)>0$, $I(x,y)>0$, and $I(x,y)=1$ or $P(x,y)=1$. If $I(x,y)=1$, then $R(x,y)=1$ and $P(x,y)=\inf\{t\in [0,1]: S(t,1)\ge 1\}=0$, but this contradicts the fact that $P(x,y)>0$. Moreover, if $P(x,y)=1$ (and $I(x,y)\neq 1$), we state that $R(x,y)=S(P(x,y),I(x,y))=1$ and, since $I(x,y)\in(0,1)$, $I(x,y)$ is a $1$-divisor. So, there is an $s_0\in (0,1)$ with $S(s_0,I(x,y))=1$. However, by the definition of $P$, we have that $1=P(x,y)=\inf\{s\in[0,1]:S(s,I(x,y)\ge 1\}\le s_0<1$, and this is a contradiction.
	\end{proof}
	\label{strongDecDrastic}
\end{proposition}

The above proposition makes us think that the relation between $0$-divisors (respectively, $1$-divisors) of a t-norm $T$ (respectively, of a t-conorm $S$) is the keystone in the existence of strong decompositions.\\
The argument in the proof is based upon avoiding   $T_D(P,I)$ being positive. Since $I$ could achieve any value (see Proposition \ref{simmin}), any of these values must have a $0$-divisor.

The next step is the generalization of Proposition \ref{strongDecDrastic} to any t-norm. For that reason, we state the following definition.

\begin{definition} Let $S$ (respectively, $T$) be a triangular conorm (respectively, a t-norm). Given any $w\in [0,1]$, we define the \emph{$1$-interval}\footnote{It is straightforward to check that $D^1_S(w)$ and $D^0_T(w)$ are intervals because of the monotonicity of $S$ and $T$. Moreover, $1\in D^1_S(w)$ and $0\in D^0_T(w)$ are always satisfied.} associated to $w$ as $D^1_S(w)=\{t\in [0,1]: S(t,w)=1\}$. Similarly, the \emph{$0$-interval} associated to $w$ is defined by $D^0_T(w)=\{t\in[0,1]: T(t,w)=0\}$.
\end{definition}

\begin{proposition} Let $S$ and $T$ respectively be a t-conorm and a t-norm, and let $X$ be a set. Every fuzzy binary relation on $X$ is strongly decomposable with respect to $S$ and $T$ if and only if for all $w\in[0,1]$ it holds true that $D^1_S(w)\cap D^0_T(w)\neq \emptyset$, and the triangular conorm $S$ is continuous on the first coordinate.
	\begin{proof} To prove the direct implication, for every $w\in[0,1]$ consider a fuzzy binary relation $R$ with $R(x,y)=1$ and $R(y,x)=w$. By hypothesis $R$ strongly decomposes into a pair $(P,I)$. We have that $1=R(x,y)=S(P(x,y),I(x,y))$ and $0=T(P(x,y),I(x,y))$, and now using that $w=I(x,y)$ (Proposition \ref{simmin}), we obtain that $P(x,y)\in D^1_S(w)\cap D^0_T(w)$.
	
	To prove the converse implication, we only need to see that the decomposition $(P,I)$ obtained in the proof of Proposition \ref{Cdecomp} satisfies $T(P(x,y),I(x,y))=0$ for all $x,y\in X$. To do so, we can assume without loss of generality that $I(x,y),P(x,y)>0$. First, when $R(x,y)=1$, we have that $P(x,y)=\inf\{t\in[0,1]: S(t,I(x,y))\ge 1\}=\inf D_S^1(I(x,y))$. In this situation, we claim that $P(x,y)=\min D_S^1(I(x,y))$ because $S$ is continuous on the first component. Moreover, since $P(x,y)$ is the minimum of $D_S^1(I(x,y))$, and $D_S^1(I(x,y))\cap D_T^0(I(x,y))\neq \emptyset$, we conclude that $P(x,y)\in D^0_T(I(x,y))$ and we get that $T(P(x,y),I(x,y))=0$. In addition, when $R(x,y)<1$, we define $\bar{R}$ as the fuzzy binary relation that coincides with $R$ everywhere except on the pair $(x,y)$, for which we define $\bar{R}(x,y)=1$. $\bar{R}$ decomposes into $(\bar{P}, I)$. It can be checked, using the definition, that $P(x,y)\le \bar{P}(x,y)$. Finally, notice that we are now in the previous situation already analyzed, because $\bar{R}(x,y)=1$, and $0=T(\bar{P}(x,y),I(x,y))\ge T(P(x,y),I(x,y))$.
	\end{proof}
	\label{ExisSD}
\end{proposition}

Taking into account the previous proof, the condition of uniqueness is quite natural. The set $D_S^1(I(x,y))\cap D_T^0(I(x,y))$ contains all the possible values for the asymmetric part $P(x,y)$. Hence, in order to achieve uniqueness, it is enough to request this set to have just a single element.

\begin{proposition} \label{45pr} Every binary fuzzy relation on $X$ is strongly decomposable with respect to $S$ and $T$, in a unique way, if and only if the t-conorm $S$ is continuous in the first component and for all $w\in[0,1]$ it holds true that $|D^1_S(w)\cap D^0_T(w)|=1$.
    \begin{proof}
        Suppose that there are $w, s, t\in [0,1]$ with $s < t$ and such that $s,t\in D_S^1(w)\cap D_T^0(w)$. We fix now a pair $x,y\in X$ and define a binary relation $R \in \mathcal{B}_X$ as $R(y,x)=w$ and $R(a,b) =  1$ if $(a,b) \neq (y,x)$. Additionally, we define three binary relations $I,P,P'\in \mathcal{B}_X$ for every $a,b\in X$ as follows: $I(a,b)=\min\{R(a,b),R(b,a)\}$, $P(a,b)=P'(a,b)=R(a,b) \searrow I(a,b)$ if $(a,b)\neq (x,y)$, $P(x,y)=s$ and $P'(x,y)=t$. It remains to prove that $R$ strongly decomposes in $(P,I)$ as well as $(P',I)$, and therefore the decomposition is not unique.\\
        Using the arguments in the proof of Proposition \ref{ExisSD} we only need to check that $R(x,y)=S(P(x,y),I(x,y))=S(P'(x,y),I(x,y))$ and $0=T(P(x,y),I(x,y))=T(P'(x,y),I(x,y))$. But these assertions are indeed true because $1=S(s,w)=S(t,w)$ and $0=T(t,w)=T(t,w)$ by hypothesis.
        
        Suppose that $R$ strongly decomposes into $(P,I)$ and $(P',I')$. Proposition \ref{simmin} guarantees that $I=I'$, so there are $x,y\in X$ such that $P(x,y)\neq P'(x,y)$. Now we define $s=P(x,y)$, $t=P'(x,y)$ and $w=I(x,y)$. To conclude, by the definition of a strong decomposition we can guarantee now that $s,t\in D^1_S(w)\cap D^0_T(w)$.
    \end{proof}
\end{proposition}

The main problem that we face when dealing with strong decompositions is that the conditions for their existence and uniqueness are too restrictive. In the \ref{app:examples}, we can see some examples for the most usual t-norms and t-conorms. For that reason, we will study now a weaker version of a decomposition defined by only using a t-conorm.\\
Furthermore, another reason to proceed in that way is that the main literature, as far as we know, uses this second type of decomposition, namely the so-called weak decomposition.

\section{Weak decompositions of fuzzy binary relations} 
\label{sec:weakDec}

We can encounter, in the specialyzed literature, some pioneer introduction of decompositions similar to the weak ones considered here, as isolated proposals to derive two binary relations from a single one generalizing the crisp decompositions. However, no general framework discussing decompositions was proposed in those studies \cite{Banerjee1994121, billot2012economic, Dutta1987215, MONTERO198645, Ovchinnikov1981169}.\\
Richardson in \cite{RICHARD2011367} proposed a general framework and it has been adopted by other authors as Fono and Andjiga in \cite{Fono2005372} or Gibilisco et al. in \cite{Gibilisco2014}. The general framework introduced by Richardson is a priori more general than Definition \ref{weakDec}. However, in order to obtain results, they request to all their decompositions to satisfy a property that they call simplicity, and then, because of having asked this extra property, our Definition \ref{weakDec} becomes more general than theirs with that condition of ``simplicity'' added.

\begin{proposition} Let $S$ be a t-conorm and $X$ a set. Every fuzzy relation $R$ on $X$ admits a  weak  decomposition if and only if $S$ is continuous in the first coordinate (with respect to the Euclidean topology on the unit interval $[0,1]$).
	\begin{proof} Notice that Proposition \ref{Cdecomp} guarantees the direct implication. For the converse implication, we also use the decomposition $(P,I)$ defined from $R$ proposed in the proof of Proposition \ref{Cdecomp}. It only remains that such pair is a weak decomposition. But this is indeed the case because, if $I(x,y)=1$ (with $x,y\in X$), then $R(x,y)=1$ and $P(x,y)=\inf\{t\in[0,1]: S(t,I(x,y))\ge R(x,y)\}=0$.
	\end{proof}
	\label{WCdecomp}
\end{proposition}

In the following example, we compute  some weak decompositions for the main t-conorms, using the formula from Proposition \ref{Cdecomp}, namely $P_S=I\searrow_S R$.

\begin{example} Let $R$ be a fuzzy binary relation on $X$. Let $S$ be a t-conorm. Using the proof of Proposition \ref{Cdecomp} we can compute a decomposition $(P_S,I)$ of $R$ with respect to several well-known t-conorms (Example \ref{example:standarnorms} contains their definitions).
	\begin{equation}
	P_{\max}(x,y)=
		\begin{cases}
	R(x,y) &\text{if } R(x,y)>R(y,x),\\
	0 & \text{otherwise}.
	\end{cases}
	\label{MaxDec}
	\end{equation}
	\begin{equation}
	P_{S_{\text{\L}}}(x,y)=\begin{cases}
	R(x,y)-R(y,x) &\text{if } R(x,y)>R(y,x),\\
	0 & \text{otherwise}.
	\end{cases}
	\label{LDec}
	\end{equation}
	\begin{equation}
	P_{S_P}(x,y)=\begin{cases}
	\frac{R(x,y)-R(y,x)}{1-R(y,x)} &\text{if } R(x,y)> R(y,x),\\
	0 &\text{otherwise}.
	\end{cases}
	\label{ProbDec}
	\end{equation}
	and $I(x,y)=\min\{R(x,y),R(y,x)\}$ (for every $x,y\in X$). We point out that these decompositions already appeared as suitable expressions in the previous literature, in this setting. For instance, $P_{S_{\max}}$ was considered in \cite{Ovchinnikov1981169,Dutta1987215}, $P_{S_{\text{\L}}}$ in \cite{Dutta1987215,MONTERO198645} and $P_{S_P}$ in \cite{Gibilisco2014}.
	\label{EDec}
\end{example}

The last proposition in this section characterizes the uniqueness of weak decompositions (provided that they exist).

\begin{proposition} Let $S$ be a t-conorm that is continuous in the first coordinate. Let $R$ be a fuzzy binary relation on $X$. $R$ has a unique weak decomposition in terms of $S$ if and only if $S$ is strictly increasing\footnote{Here, being strict increasing has to be interpreted in the sense of Definition \ref{DefStrict}. No t-conorm is strictly increasing on the whole domain $[0,1]^2$.} in the first coordinate.

	\begin{proof} Suppose that $R$ decomposes into two distinct weak decompositions $(P,I)$ and $(P',I)$. We can suppose without loss of generality that there exist $x,y\in X$ such that $P(x,y)>P'(x,y)$. Also $I(x,y)\neq 1$ because $P(x,y)>0$. Hence, using the condition of strict increasingness, we get that $R(x,y)=S(P(x,y),I(x,y))>S(P'(x,y),I(x,y))=R(x,y)$. Thus we have obtained a contradiction, and, consequently, the decomposition is unique.
	
	If $S$ is not strict, there are $s,t,w\in (0,1)$ with $S(s,w)=S(t,w)$ and $s\neq t$. We can define a binary relation $R$ being $R(x,y)=S(s,w)$, $R(y,x)=w$ for some $x,y\in X$ and $R(a,b) = 1$ at any other pair $(a,b)$. It can be checked that $R$ weakly decompose into $(P,I)$ and $(P',I)$ defined as $I(a,b)=\min\{R(a,b),R(b,a)\}$ (for every $a,b\in X$) and $P(x,y)=s$, whereas $P'(x,y)=t$ and  $P'(a,b) = 0$ at any other pair $(a,b)$.
	\end{proof}
	\label{UDesc}
\end{proposition}

In the Appendix, Table \ref{TableDecompositions} furnishes information about decomposability, analyzing which ones among the main t-conorms can decompose all fuzzy binary relations.

\section{Applications of   decompositions of fuzzy binary relations into preference modelling under uncertainty} \label{sec:decompositionPreferences}

In this section, we will apply the results of the previous sections about decomposition of fuzzy binary relations to fuzzy Arrovian models.

In this literature, it is required to obtain a strict preference $P$ and an indifference preference $I$ from a weak preference $R$. However, the inclination is defining a rule to obtain $P$ and $I$ firstly and, secondly, studying their properties \cite{Banerjee1994121, Gibilisco2014, Ovchinnikov1981169}. Moreover, most of the authors in Arrovian literature have used the type of decompositions studied in Sections \ref{sec:strongDec} and \ref{sec:weakDec} to define such rules \cite{Dutta1987215, Fono2005372, Gibilisco2014, Richardson1998359}. For instance, the probabilistic t-conorm induces a unique decomposition (see Equation \ref{ProbDec} in Example \ref{EDec} and Table \ref{TableDecompositions}). In the following pages, we will see that the decomposition generated by the probabilistic t-conorm satisfies the desired properties in fuzzy Social Choice. However, other conorms behave differently. Consider the following example:

\begin{example} Consider de following t-conom $S$ defined as the ordinal sum of a \L ukasiewicz t-conorm (see \cite[Definition 3.68]{Beliakov2007}).
\begin{equation*}
    S(x,y)=\begin{cases}
        \min\{0.5,x+y\} &\text{if } x,y\le 0.5,\\
        \max\{x,y\} & \text{otherwise.}
    \end{cases}
\end{equation*}
Any fuzzy binary relation weakly decomposes with respect $S$ because it is continuous (see Proposition \ref{WCdecomp}). However, consider the binary relation $R$ on $X=\{x,y\}$ defined as $R(x,y)=R(y,x)=0.5$ and $R(x,x)=R(y,y)=1$. It decomposes into $(P,I)$ where $I(x,x)=I(y,y)=1$, $P(x,y)=I(x,y)=0.5$ and $P(y,x)=0$. 
How could it be that, despite $x$ and $y$ playing the same role in $R$, $x$ is strictly preferred over $y$? This is one of the situations that are avoided implicitly in fuzzy Arrovian models \cite{Gibilisco2014}.

\end{example}

Therefore, we will first define what we expect (as Social Choice practitioners) from $R$, $P$, $I$ and the relationship among them and, straightaway, how strong and weak decompositions can be used to obtain suitable decomposition rules.

With that purpose in mind, we use the concept of preference structure (see e.g. \cite{preferenceRoubens, VanDeWalle1998105}) because its definition does not rely upon any concept related to decomposition rules. A preference structure is a triplet $(R,P,I)$ of fuzzy binary relations satisfying the properties of Definition \ref{DFpreference} below.

\begin{definition}[{{\cite[Definition 2.15]{arrowfuzzyA1}}}] Given a set $X$, a \emph{fuzzy preference} on $X$ is a triplet $(R,P,I)$ of fuzzy binary relations on $X$ that satisfies the following conditions:
	\begin{itemize}
		\item [(FP1)] $P(x,y)>0$ implies $P(y,x)=0$, for all $x,y\in X$ ($P$ is \emph{asymmetric}),
		\item [(FP2)] $I(x,y)=I(y,x)$, for all $x,y\in X$ ($I$ is \emph{symmetric}),
		\item [(FP3)] $P(x,y)\le R(x,y)$ for every $x,y\in X$,
		\item [(FP4)] $R(x,y)>R(y,x)$ if and only if $P(x,y)>0$, for every $x,y\in X$,
		\item [(FP5)] if $P(x,y)=0$ then $R(x,y)=I(x,y)$, for all $x,y\in X$,
		\item [(FP6)] if $I(x,y) \le I(z,w)$ and $P(x,y)\le P(z,w)$ then $R(x,y)\le R(z,w)$ holds true for every $x,y,z,w\in X$.
	\end{itemize}
	We call to $R$ the \emph{weak preference} component, to $P$ the \emph{strict preference} component and to $I$ the \emph{indifference} component.
	\label{DFpreference}
\end{definition}

Using this definition, we state a framework to study the decomposition rules. It is equivalent to the standard definitions in fuzzy Arrovian models \cite{Dutta1987215, Fono2005372, Gibilisco2014}. Now, we can use this framework to apply the decompositions studied in the previous sections to decomposition rules. First, we need to define what we understand as a decomposition rule:

\begin{definition} Let $\mathcal{B}$ be a family of fuzzy binary relations on a set $X$. A \emph{decomposition rule} on $\mathcal{B}$ is a map $\phi: \mathcal{B} \rightarrow \mathcal{BR}_X \times \mathcal{BR}_X$ such that for every $R \in \mathcal{B}$, if $\phi(R)=(P_R,I_R)$, then $\Lambda_R = (R,P_R,I_R)$ is a fuzzy preference in the sense of Definition \ref{DFpreference}.  
	\label{DescD}
\end{definition}

Our goal in this section will be to study when the decompositions already introduced through sections \ref{sec:strongDec} and \ref{sec:weakDec} induce or are compatible with a decomposition rule, in the following sense of the next Definition \ref{piol}.

\begin{definition} \label{piol} Let $T$ be a t-norm, $S$ a t-conorm and $\phi$ a decomposition rule on $\mathcal{B}\subseteq \mathcal{BR}_X$.
\begin{itemize}
    \item [(i)] We say that $\phi$ is \emph{$(T,S)$-strong (resp. $S$-weak) compatible} if, for every $R\in \mathcal{B}$, it holds true that $\phi(R)$ is a strong (resp. weak) decomposition of $R$ with respect to $T$ and $S$ (resp. $S$).
    \item [(ii)] We say that $\phi$ is \emph{induced by $(T,S)$ (resp $S$)} if $\phi$ is $(T,S)$-strong (resp. $S$-weak) compatible and for every $R\in \mathcal{B}$ it holds, in addition, that $\phi(R)$ is the unique decomposition of $R$ which induces a preference. That is, if $(P,I)$ is a decomposition of $R$ and $(R,P,I)$ is a preference, then $\phi(R)=(P,I)$. In that case, we denote by $\phi_{T,S}$ (resp. $\phi_S$) the induced decomposition rule.
\end{itemize}
\end{definition}

\begin{example} The decomposition rule $\phi_1(R)=(P_R^1,I_R)$ defined in $\mathcal{B}=\mathcal{BR}_X$ as $I_R(x,y)=\min\{R(x,y),R(y,x)\}$ and $P_R^1(x,y)=\min\{0,\frac{1}{2}(R(x,y)-R(y,x))\}$ is neither $(T,S)$-strong nor $S$-weak compatible for any $T$ and $S$: if $R(x,y)=1$ and $R(y,x)=0$ then, $R(x,y)\neq S(P_R^1(x,y),I_R(x,y))$.

Moreover, consider the decompositions rules $\phi_2=({P_{S_{\text{\L}}}},I_R)$ (see Example \ref{EDec}) and $\phi_3=(P'_{S_{\text{\L}}}, I_R)$ defined as ${P'_{S_{\text{\L}}}}_R(x,y)=1$ if $1=R(x,y)>R(y,x)$, and $P'_{S_{\text{\L}}}(x,y)=P_{S_{\text{\L}}}(x,y)$ otherwise. They are both $S_{\text{\L}}$-weak compatible, so neither of them are induced by $S_{\text{\L}}$.
Finally, Proposition \ref{UMax} proves that $\phi_{S_{\max}}$ is induced by $S_{\max}$. 

\end{example}

First, we will focus on the case $(i)$ in the Definition \ref{piol} above. We will study which t-norms and t-conorms admit decomposition rules compatible with them. Next, we will focus on the t-norms and t-conorms which induce decomposition rules (case $(ii)$).

As in the previous sections, we will restrict our study to the case $\mathcal{B}=\mathcal{BR}_X$, since other domains could lead to different results.

\begin{proposition} Let $R$ be a fuzzy binary relation on $X$. Let $S$ be a t-conorm and let $T$ be a t-norm. If $(P,I)$ is a strong or a weak decomposition of $R$ with respect to $T$ and $S$ or respect to $S$, then the  properties $FP1$, $FP2$, $FP3$, $FP5$, $FP6$ and the direct implication   in the statement $FP4$ of Definition \ref{DFpreference} hold true.\\
Moreover, the converse of $FP4$ is satisfied if for every $t\in [0,1)$ there exist a neighborhood of $0$ in which $S(\cdot,t)$ is strictly increasing.
	\begin{proof} Conditions $FP1$ and $FP2$ are obtained from the definition of a decomposition. To prove $FP3$ notice that $R(x,y)=S(P(x,y), I(x,y)) \ge S(P(x,y), 0)\allowbreak=P(x,y)$. 
		
	To prove $FP5$ notice that if $P(x,y)=0$, then $R(x,y)=S( P(x,y),I(x,y))\allowbreak =S(0, I(x,y))=I(x,y)$.
		
	Property $FP6$ is a direct consequence of the monotonicity of $S$.
		
	Finally, to prove the direct implication of $FP4$, observe that if $R(x,y)>R(y,x)$, by property $FP6$ it follows now that $[I(x,y)>I(y,x)] \vee [P(x,y)>P(y,x)]$ holds true.  But, since $I$ is symmetric, so that $I(x,y)=I(y,x)$, we get $P(x,y)>P(y,x)\ge 0$.
		
	Additionally, suppose that $S$ satisfies the last property. If $P(x,y)>0$, then $P(y,x)=0$. Hence we  may now conclude that $R(y,x)=S(0,I(y,x))<S(P(x,y),I(x,y))=R(x,y)$.
	\end{proof}
	\label{PropDesc}
\end{proposition}

In the following example, we can see a decomposition that satisfies all the properties above but the converse implication in $FP4$ and, consequently, does not give rise to a fuzzy preference.

\begin{example} Let $R$ be a fuzzy binary relation defined on a set $X$. Consider the fuzzy binary relations $P$ and $I$ on $X$ defined as:
	\begin{equation*}
		P(x,y)=\begin{cases} 
			R(x,y) & \text{if } R(x,y)>R(y,x),\\
			R(x,y) & \text{if } 1>R(x,y)= R(y,x),\\
			0 &\text{otherwise}.
		\end{cases}
	\end{equation*}
and $I(x,y)=\min\{R(x,y),R(y,x)\}$. Notice that $R$ weakly decomposes into $(P,I)$ with respect $S_{\max}$ and this decomposition is not the same as the one arising in Example \ref{EDec}.
\label{ex:twodecompositionrules} 
\end{example}

The case of a weak decomposition with respect $S_{\max}$ is paradigmatic in the literature, and we will use it to motivate the main result in this section.\\

The decomposition \ref{MaxDec} shown in Example \ref{EDec} actually appeared early in the specialized literature \cite{Ovchinnikov1981169}, usually as a suitable decomposition rule in several contexts. The example above shows that it is not the unique weak decomposition with respect to $S_{\max}$. However, the following proposition proves that it is the unique decomposition that defines a fuzzy preference. That is, $S_{\max}$ induces a decomposition rule. 

\begin{proposition} Let $R$ be a fuzzy binary relation on a set $X$. There exists a unique decomposition $(P,I)$ of $R$ with respect to the t-conorm $S_{\max}$ satisfying that $(R,P,I)$ is a fuzzy preference.
	\begin{proof} Consider a fuzzy binary relation $R$. The maximum is a continuous t-conorm. By   Proposition \ref{Cdecomp}, $R$ is decomposable. Suppose that there are two different decompositions $(P,I)$ and $(P',I')$ of $R$ such that both $(R,P,I)$ and $(R,P',I')$ are fuzzy preferences. First, by   Proposition \ref{simmin}, $I=I'$ holds. If $P\neq P'$, we can assume, without loss of generality, that there are $x,y\in X$ with $P(x,y)>P'(x,y)$. \\
	From the equality $R(x,y)=\max\{P(x,y),I(x,y)\}=\max\{P'(x,y),I(x,y)\}$ we obtain that $R(x,y)=I(x,y)$. Using that $I(x,y)=\min\{R(x,y),R(y,x)\}$ (Proposition \ref{simmin}), we conclude that $R(x,y)\le R(y,x)$. We have finally arrived at a contradiction because, using $FP4$, we obtain that $P(x,y)=0$.
	\end{proof}
	\label{UMax}
\end{proposition}

\begin{remark} In Proposition \ref{PropDesc} we have proved that the property of $S$ \emph{``for every $t$, having a neighbourhood of $0$ such that $S(t,\cdot)$ is strictly increasing''} implies $FP4$.

If we could find a property of $S$ equivalent to $FP4$, then the characterization of decompositions that are a preference would be completed. Unfortunately, Example \ref{ex:twodecompositionrules} makes us think that such property does not exist.

Assuming the existence of such property $q$, $S_{\max}$ may or may not satisfy $q$. If $S_{\max}$ satisfied it, then all weak decompositions (with respect $S_{\max}$) should define a preference, but the decomposition in Example \ref{ex:twodecompositionrules} does not define a preference. Conversely, in case $S_{\max}$ satisfied $q$, then no weak decomposition should define a preference, but, according with Proposition \ref{UMax}, this is not true.

Then, such property $q$ should not exist. We conclude that an equivalence to $FP4$ relies on something more than some properties of t-conorms. 
\end{remark}


The following proposition generalizes the previous example to any t-conorm that may be used in weak decompositions.

\begin{proposition}\label{mj} Let $S$ be a right-continuous t-conorm.
$S$ induces a weak decomposition rule if, and only if, for every $w,t,s\in [0,1]$ with $t\neq s$ it holds true that,  if $S(t,w)=S(s,w)$, then $S(t,w)=w$.
\begin{proof} To prove the direct implication, suppose that there are $w,s,t\in [0,1]$ such that $S(t,w)=S(s,w)>w$, and we will prove that there is an $R\in \mathcal{BR}_X$ which has at least two decompositions which induce fuzzy preferences.\\
To see that, given a pair $a,b\in X$ define a fuzzy binary relation $R$ satisfying $R(a,b)=S(t,w)$;    $R(b,a)=w$ and $R(c,d)=1$ if $(c,d) \neq (a,b)$ and also $(c,d) \neq (b,a)$. Consider now the three fuzzy binary relations $I(x,y)=\min\{R(x,y),R(y,x)\}$ for all $x,y\in X$, $P(x,y)=P'(x,y)=0$ if $(x,y)\neq(a,b)$, $P(a,b)=t$ and $P'(a,b)=s$.\\
Now we can assert that $(P,I)$ and $(P',I)$ are weak decompositions of $R$. Moreover $(R,P,I)$ and $(R,P',I)$ are preferences because they are defined by means of a decomposition and it is immediate to check that they indeed satisfy $FP4$. Therefore, $S$ can not induce any decomposition rule.

Concerning the converse implication, we may notice that the unique property of the maximum used in the proof of Proposition \ref{UMax} is just our hypothesis here. We could transcript the same proof after interchanging \emph{"From the equality $R(x,y)=\max\{P(x,y),I(x,y)\}=\max\{P'(x,y),I(x,y)\}$ we obtain that $R(x,y)=I(x,y)$"} with \emph{"From the equality $R(x,y)=S(P(x,y),I(x,y))=S(P'(x,y),I(x,y))$ we obtain that $R(x,y)=I(x,y)$"}.
\end{proof}
\end{proposition}

\begin{remark} Notice that if $S$ is strictly increasing in the first coordinate, $S$ satisfies the hypothesis of the Proposition \ref{mj} because if $S(t,w)=S(s,w)$ then $w=1$. This fact is consistent with the combination of Propositions \ref{UDesc} and \ref{PropDesc}, which also guarantees the existence of the induced decomposition rule $\phi_S$.

\end{remark}


\section{Conclusion}
The concept of a fuzzy preference, defined as a triplet that reminds us the usual weak preference, strict preference and indifference arising in the crisp setting, has been analyzed from the point of view of the existence of decompositions in which the fuzzy weak preference generates, in a way, the associated strict fuzzy preference as well as the fuzzy indifference. Whereas in the crisp setting the decompositions are always unique, this is no longer true in the fuzzy setting. Consequently, we have also analyzed questions related to the existence and uniqueness of decompositions of fuzzy preferences, as shown in the main sections of the present manuscript.

\section*{Acknowledgments}


This work has been partially supported by the  research project PID2019-108392GB-I00 (AEI/10.13039/ $ $501100011033), the \emph{Ayudas para la Recualificación del Sistema Universitario Español para 2021-2023, UPNA. Modalidad Margarita Salas} funded by the European Union - NextGenerationEU, and a predoctoral grant from the UPNA Research Institutes.

\bibliographystyle{abbrv} 
\bibliography{mybibfile}

\appendix

\section{Decompositions in restricted domains}
\label{sec:restrictDomains}

As we have discussed in the Introduction, in most  problems related to preferences, we require $R$ to satisfy additional properties. Usually, they are some kind of transitivity and connectednes.\\
So, instead of studying the existence of decompositions on $\mathcal{BR}_X$, we could study them on the subsets defined by those transitivities and connectedness assumptions.   In this section, we will give some geometric ideas about the conditions of existence for the decomposition in such domains, but we will not furnish any formal proof.\\

Consider ${\mathcal{BR}_X}^{T'}$ (resp. ${\mathcal{BR}_X}^{S'}$) the set of all $T'$-transitive (resp. $S'$-connected)\footnote{A fuzzy binary relation $R$ is said $T'$-transitive if for every $x,y,z\in X$ the inequality $R(x,z)\ge T'(R(x,y),R(y,z)$ is satified. Morover, it is said $S'$-connected if for every $x,y\in X$ $S'(R(x,y),R(y,x))=1$. See \cite[Chapter 3]{Gibilisco2014} or \cite{Raventos-Pujol202003} for a more extended exposition.} binary relations of $\mathcal{BR}_X$. Notice that, a priori, the t-norm and the t-conorm defining the transitivity and the connectedness and the ones defining the decomposition may or may not coincide. For that reason, we use $T'$ and $S'$ to denote the ones corresponding to the transitivity and the connectedness, and $T$ and $S$ the ones defining the decomposition.\\

First, if $\phi$ is a decomposition rule on $\mathcal{BR}_X$, then its restriction to one of these sets will also be a decomposition rule.\\
However, could it be the case that, for example, there is no decomposition rule $\phi$ $(T,S)$-strong compatible on $\mathcal{BR}_X$ but there are such rules\footnote{We can make analogous suggestions about $S$-compatibility and decomposition rules induced by $(T,S)$ or $S$.} on ${\mathcal{BR}_X}^{T'}$ or on ${\mathcal{BR}_X}^{S'}$?\\

\begin{figure}[ht]
    \centering
    \begin{subfigure}[t]{0.3\textwidth}
        \centering
        \begin{tikzpicture}[scale=0.5, every node/.style={scale=1}]

\begin{axis}[axis lines=middle,
            xlabel={$R(x,y)$},
            ylabel={$R(y,x)$},
            enlargelimits,
            ytick={1},
            xtick={1},
            yticklabels={$1$},
            xticklabels={$1$},]
            
\addplot[name path=F,red,domain={0:1}] {1} node[pos=.8, above]{};
\addplot[name path=F,red,domain={0:1}] {0} node[pos=.8, above]{};

\addplot [mark=none, red] coordinates {(0, 0) (0, 1)};
\addplot [mark=none, red] coordinates {(1, 0) (1, 1)};

\end{axis}

\node[red] at (3,3) {$D_{S_D}$};

\end{tikzpicture}
    \end{subfigure}
    \begin{subfigure}[t]{0.3\textwidth}
        \centering
        \begin{tikzpicture}[scale=0.5, every node/.style={scale=1}]

\begin{axis}[axis lines=middle,
            xlabel={$R(x,y)$},
            ylabel={$R(y,x)$},
            enlargelimits,
            ytick={1},
            xtick={1},
            yticklabels={$1$},
            xticklabels={$1$},]

\addplot[name path=F,red,domain={0:1}] {x} node[pos=.8, above]{};
\addplot[name path=F,red,domain={0:1}] {0} node[pos=.8, above]{};

\addplot [mark=none, red] coordinates {(0, 0) (0, 1)};

\end{axis}

\node[red] at (3,4) {$D_{T_P,S_{\text{\L}}}$};

\end{tikzpicture}
    \end{subfigure}
    \begin{subfigure}[t]{0.3\textwidth}
        \centering
        \begin{tikzpicture}[scale=0.5, every node/.style={scale=1}]

\begin{axis}[axis lines=middle,
            xlabel={\tiny $R(x,y)$},
            ylabel={\tiny $R(y,x)$},
            enlargelimits,
            ytick={1,0.5, 0.75},
            xtick={1,0.5, 0.75},
            yticklabels={\tiny $1$, \tiny $\frac{1}{2}$, \tiny $\frac{3}{4}$},
            xticklabels={\tiny $1$, \tiny $\frac{1}{2}$, \tiny $\frac{3}{4}$},]

\addplot[name path=F,red,domain={0:1}] {x^2-x+1} node[pos=.8, above]{};

\addplot[name path=G1,red,domain={0.75:1}] {(1+sqrt(4*x-3))*0.5}node[pos=.1, below]{};
\addplot[name path=G2,red,domain={0.75:1}] {(1-sqrt(4*x-3))*0.5}node[pos=.1, below]{};

\addplot[name path=B, opacity = 0, domain={0:1}] {0}node[pos=.1, below]{$g$};

\addplot[color=red!10, opacity = 0.4]fill between[of=F and G1, soft clip={domain=0.75:1}];

\addplot[color=red!10, opacity = 0.4]fill between[of=F and B, soft clip={domain=0:0.75}];

\addplot[color=red!10, opacity = 0.4]fill between[of=G2 and B, soft clip={domain=0.75:1}];

\end{axis}

\node[red] at (3,3) {$D_{T_{\text{\L}}, S_P}$};

\end{tikzpicture}
    \end{subfigure}
    \caption{Representations of $D_{S_D}$, $D_{T_P,S_{\text{\L}}}$ and $D_{T_{\text{\L}},S_P}$. Since none of these sets is the unit square $[0,1]^2$, these combinations of t-norms and t-conorms do not allow decompositions (compare with Table \ref{TableDecompositions}).}
    \label{Cuadro1}
\end{figure}
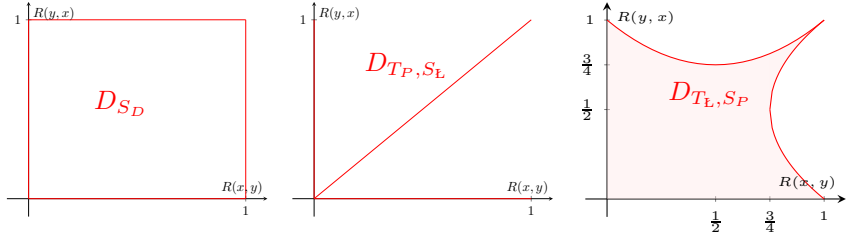

Before we go further  on the possible answer to these questions, we need to introduce some terminology. First, as we have noticed before, the existence of decompositions and their uniqueness of a binary relation $R$ depend only on the comparison between the values of $R(x,y)$ and $R(y,x)$ for every $x,y\in X$. We can illustrate this idea using the sets $D_S$ and $D_{T,S}$ defined below (see some examples in Figure \ref{Cuadro1}):

The points $a,b\in [0,1]^2$ are the combinations of the values for $R(x,y)$ and $R(y,x)$ in which a decomposition is possible. So, we can state that a binary relation $R$ is strong decomposable with respect $T$ and $S$ (resp. weak decomposable with respect $S$) if and only if for every $x,y\in X$ $(R(x,y),R(y,x))\in D_{T,S}$ (resp. $D_S$). Moreover, since $ \mathcal{BR}_X$ imposes no restrictions over the values, we can conclude that all binary relations are strongly (resp. weakly) decomposable if $D_{T,S}=[0,1]^2$ (resp. if $D_{S}=[0,1]^2$).\\

Now, we can state that in $\mathcal{BR}_X^{T'}$ apply the same results we have obtained for $\mathcal{BR}_X$. First, $T'$-transitivity does not impose restrictions on the possible relations between $R(x,y)$ and $R(y,x)$. In other words, given $(a,b)\in [0,1]^2$, there is a $T'$-transitive binary relation $R$ with $R(x,y)=a$ and $R(y,x)=b$. So, we can conclude that all binary relations in $\mathcal{BR}_X$ decompose if, and only if, all binary relations in ${\mathcal{BR}_X}^{T'}$ decompose.\\

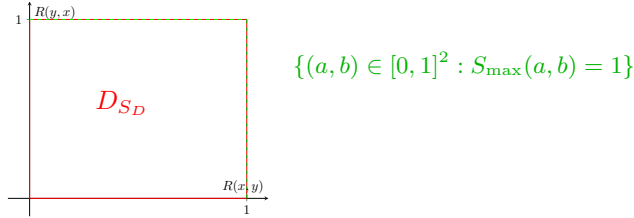
\begin{figure}[ht]
    \centering
    \begin{tikzpicture}[scale=0.5, every node/.style={scale=1}]

\begin{axis}[axis lines=middle,
            xlabel={$R(x,y)$},
            ylabel={$R(y,x)$},
            enlargelimits,
            ytick={1},
            xtick={1},
            yticklabels={$1$},
            xticklabels={$1$},]

\addplot[name path=F,red,domain={0:1}] {1} node[pos=.8, above]{};
\addplot[name path=F,red,domain={0:1}] {0} node[pos=.8, above]{};

\addplot [mark=none, red] coordinates {(0, 0) (0, 1)};
\addplot [mark=none, red] coordinates {(1, 0) (1, 1)};

\addplot[name path=F,green, dashed ,domain={0:1}] {1} node[pos=.8, above]{};
\addplot [mark=none, black!30!green, dashed] coordinates {(1, 0) (1, 1)};

\end{axis}

\node[red] at (3,3) {$D_{S_D}$};
\node[black!30!green] at (12,4) {\small$\{(a,b)\in [0,1]^2: S_{\max}(a,b)=1$\}};

\end{tikzpicture}
    \caption{All binary relations in ${\mathcal{BR}_X}^{S_{\max}}$ weakly decompose with respect $S_D$ because $\{(a,b)\in [0,1]^2: S_{\max}(a,b)=1\}\subseteq D_{S_D}$. However, not every binary relation in $\mathcal{BR}_X$ weakly decompose with respect $S_D$.}
    \label{Cuadro2}
\end{figure}

However, the restrictions imposed by connectedness are different. In this case, $S'$-completeness imposes that $S'(R(x,y),R(y,x))=1$. Since now some combinations of values will not be feasible, we conclude that all binary relations on ${\mathcal{BR}_X}^{S'}$ will decompose if, and only if, $\{(a,b)\in [0,1]^2: S'(a,b)=1\} \subseteq D_S$ (resp. if, and only if, $\{(a,b)\in [0,1]^2: S'(a,b)=1\} \subseteq D_{T,S}$). We can see an example in Figure \ref{Cuadro2}.

\section{Some examples using the main t-norms and t-conorms}
\label{app:examples}

In this section, we have included the definitions of the most used t-norms and t-conorms \cite{Beliakov2007}. We have also included in some tables which combinations allow decompositions and which ones are compatible and induce decomposition rules.

\begin{example} \label{example:standarnorms} Here we give a brief account about the most widespread t-norms and their respective dual t-conorms:\\
	\begin{itemize}
		\item [] The minimum $T_{\min}$ and the maximum $S_{\max}$:
		\begin{equation*}
		T_{\min}(x,y)=\min\{x,y\}, \qquad S_{\max}(x,y)=\max\{x,y\}.
		\end{equation*}
		\item [] The product $T_P$ and the probabilistic sum $S_P$:
		\begin{equation*}
		T_P(x,y)=xy, \qquad S_P(x,y)=x+y-xy.
		\end{equation*}
		\item [] The \L ukasiewicz t-norm $T_{\text{\L}}$ and the \L ukasiewicz t-conorm $S_{\text{\L}}$:
				\begin{equation*}
				T_{\text{\L}}(x,y)= \max\{0,x+y-1\}, \qquad S_{\text{\L}}(x,y)= \min\{1,x+y\}.
				\end{equation*}
		\item [] The drastic product $T_D$ and drastic sum $S_D$:
		\begin{equation*}
			T_D(x,y)=\begin{cases}
					x &\text{if } y=1,\\
					y &\text{if } x=1,\\
					0 &\text{otherwise.}
				\end{cases} \qquad \qquad
			S_D(x,y)=\begin{cases}
					x &\text{if } y=0,\\
					y &\text{if } x=0,\\
					1 &\text{otherwise.}
				\end{cases}
		\end{equation*}
		\end{itemize}

And the associated divisor intervals are:

\begin{equation*}
D^1_{S_D}(w)=\begin{cases}
\{1\} &\text{if } w=0,\\
(0,1] &\text{if } w\notin \{0,1\},\\
[0,1] &\text{if } w=1.
\end{cases} \qquad \qquad
D^0_{T_D}(w)=\begin{cases}
\{0\} &\text{if } w=1,\\
[0,1) &\text{if } w\notin \{0,1\},\\
[0,1] &\text{if } w=0.
\end{cases}
\end{equation*}
\begin{equation*}
D^1_{\max}(w)=D^1_{S_P}(w)=\begin{cases}
\{1\} &\text{if } w\neq 1,\\
[0,1] &\text{if } w=1.
\end{cases} \qquad \qquad
D^0_{\min}(w)=D^0_{T_P}(w)=\begin{cases}
\{0\} &\text{if } w\neq 0,\\
[0,1] &\text{if } w=0. 
\end{cases}
\end{equation*}
\begin{equation*}
D^1_{S_{\text{\L}}}(w)=[1-w,1], \qquad \qquad D^0_{T_{\text{\L}}}(w)=[0,1-w].
\end{equation*}
\end{example}

\begin{example} The family $(T_\lambda^{SS})_{\lambda\in [-\infty,\infty]}$ of Schweizer-Sklar t-norms is defined as:
\begin{equation*}
    T_\lambda^{SS}(x,y)=\begin{cases}
        \min\{x,y\} & \text{if } \lambda = -\infty, \\
        T_P(x,y) & \text{if } \lambda =0, \\
        T_D(x,y) & \text{if } \lambda = \infty,\\
        (\max\{x^\lambda+y^\lambda -1,0\})^{\frac{1}{\lambda}} & \text{otherwise.}
    \end{cases}
\end{equation*}
And the respective family $(S_\lambda^{SS})_{\lambda\in [-\infty,\infty]}$ of t-conorms is defined as:
\begin{equation*}
    S_\lambda^{SS}(x,y)=\begin{cases}
        \max\{x,y\} & \text{if } \lambda = -\infty, \\
        S_P(x,y) & \text{if } \lambda =0, \\
        S_D(x,y) & \text{if } \lambda = \infty,\\
        1-\left(\max\{(1-x)^\lambda+(1-y)^\lambda -1,0\}\right)^{\frac{1}{\lambda}} & \text{otherwise.}
    \end{cases}
\end{equation*}

For $\lambda\neq -\infty,0,\infty$ we have that:
\begin{equation*}
    D_{T_\lambda^{SS}}^0(w)=\begin{cases}
        \{0\} &\text{if } \lambda<0 \text{ and } w\neq0,\\
        [0,1] &\text{if } \lambda<0 \text{ and } w = 0, \\
        [0,(1-w^\lambda)^\frac{1}{\lambda}] &\text{otherwise}.
    \end{cases} 
\end{equation*}
\begin{equation*}
    D_{S_\lambda^{SS}}^1(w)=\begin{cases}
        \{1\} &\text{if } \lambda<0 \text{ and } w\neq1,\\
        [0,1] &\text{if } \lambda<0 \text{ and } w = 1, \\
        [1-(1-(1-w)^\lambda)^\frac{1}{\lambda},1] &\text{otherwise}.
    \end{cases}        
\end{equation*}
\end{example}

\begin{example} The family $(T_\lambda^{H})_{\lambda\in [0,\infty]}$ of Hamacher t-norms is defined as:
\begin{equation*}
    T_\lambda^{H}(x,y)=\begin{cases}
        T_D(x,y) & \text{if } \lambda = \infty, \\
        0 & \text{if } \lambda = x= y =0, \\
        \frac{xy}{\lambda+(1-\lambda)(x+y-xy)} & \text{otherwise.}
    \end{cases}
\end{equation*}

And the respective family $(S_\lambda^{H})_{\lambda\in [0,\infty]}$ of t-conorms is defined as:

\begin{equation*}
    S_\lambda^{H}(x,y)=\begin{cases}
        S_D(x,y) & \text{if } \lambda = \infty, \\
        1 & \text{if } \lambda = x= y =0, \\
        \frac{x+y-xy-(1-\lambda)xy}{1-(1-\lambda)xy} & \text{otherwise.}
    \end{cases}
\end{equation*}

If $\lambda \neq \infty$ we have that $D_{T_\lambda^{H}}^0(w)=D_{\min}^0(w)$ and $D_{S_\lambda^{H}}^1(w)=D_{\max}^1(w)$.
\end{example}

We have depicted the existence and uniqueness of strong and weak decompositions in $\mathcal{BR}_X$ for these t-norms and t-conorms in Table \ref{TableDecompositions}. For that purpose, we have applied Propositions \ref{ExisSD}, \ref{45pr}, \ref{WCdecomp} and \ref{UDesc}.\\

In Table \ref{TablePreferences} we have applied Propositions \ref{PropDesc}, \ref{UMax} and \ref{mj} to the same t-norms and t-conorms obtaining the corresponding results for decomposition rules.

\begin{landscape}

    \centering
    \begin{table}[]
        \begin{tabular}{l|c|c|c|c|c|c|}
         \diagbox[innerwidth=2.2cm]{$T$}{$S$} & Drastic & Maximum & \L ukasiewicz & Probabilistic & Schweizer-Sklar & Hamacher \\ \hline
         Drastic & $\nexists$ & $\nexists$ & $\exists$ & $\nexists$ & \makecell{$\lambda\le 0 \Rightarrow \nexists$ \\ $0<\lambda<+\infty \Rightarrow \exists$ \\ $\lambda = +\infty \Rightarrow \nexists$ } & $\nexists$ \\ \hline
         Minimum & $\nexists$ & $\nexists$ & $\nexists$ & $\nexists$ & $\nexists$ & $\nexists$\\ \hline
         \L ukasievicz & $\nexists$ & $\nexists$ & $\exists !$ & $\nexists$ & \makecell{ $\lambda \le 0 \Rightarrow \nexists$ \\ $0<\lambda <1 \Rightarrow \exists$ \\ $\lambda = 1 \Rightarrow \exists!$ \\ $\lambda >1 \Rightarrow \nexists$} & $\nexists$ \\ \hline
         Product & $\nexists$ & $\nexists$ & $\nexists$ & $\nexists$ & $\nexists$ & $\nexists$ \\ \hline
         Schweizer-Sklar & $\nexists$ & $\nexists$ & \makecell{ $\lambda<1 \Rightarrow \nexists$ \\ $\lambda=1 \Rightarrow \exists!$ \\ $\lambda >1 \Rightarrow \exists$}& $\nexists$ & \makecell{$\lambda<1 \Rightarrow \nexists$ \\ $\lambda = 1 \Rightarrow \exists !$ \\ $\lambda > 1 \Rightarrow \exists$ \\ $\lambda=+\infty \Rightarrow \nexists$} & $\nexists$ \\ \hline
         Hamacher & $\nexists$ & $\nexists$ & $\nexists$ & $\nexists$ & $\nexists$ & $\nexists$ \\
         
         \hline \hline
         \begin{tabular}{@{}l@{}} Weak \\ decomposition \end{tabular} & $\nexists$ & $\exists$ & $\exists$ & $\exists !$ & \makecell{$\lambda=-\infty \Rightarrow \exists$ \\ $-\infty < \lambda \le 0 \Rightarrow \exists !$ \\ $0<\lambda<+\infty \Rightarrow \exists$ \\ $\lambda = +\infty \Rightarrow \nexists$} & \makecell{$\lambda < +\infty \Rightarrow \exists !$ \\ $\lambda = +\infty \Rightarrow \nexists$} \\
         \hline
        \end{tabular}
        \caption{Decompositions for main t-norms and t-conorms. The symbol $\nexists$ means that the correspondig decomposition does not exist, $\exists$ means that it exist but it is not unique and $\exists !$ means that it exists a unique decomposition. Moreover, the characterization of  the existence and the uniqueness of the strong decomposition of Schweizer-Skalar t-norms and t-conorms corresponds to the combination $(T_{\lambda}^{SS},S_{\lambda}^{SS})$ with the same parameter $\lambda$.}
        \label{TableDecompositions}
    \end{table}

    \begin{table}[]
        \begin{tabular}{l|c|c|c|c|c|c|}
         \diagbox[innerwidth=2.2cm]{$T$}{$S$} & Drastic & Maximum & \L ukasiewicz & Probabilistic & Schweizer-Sklar & Hamacher \\ \hline
         Drastic & $\nexists$ & $\nexists$ & $\exists CR$ & $\nexists$ & \makecell{$\lambda\le 0 \Rightarrow \nexists$ \\ $0<\lambda<+\infty \Rightarrow$ ? \\ $\lambda = +\infty \Rightarrow \nexists$ } & $\nexists$ \\ \hline
         Minimum & $\nexists$ & $\nexists$ & $\nexists$ & $\nexists$ & $\nexists$ & $\nexists$ \\ \hline
         \L ukasievicz & $\nexists$ & $\nexists$ & $IDR$ & $\nexists$ & \makecell{ $\lambda \le 0 \Rightarrow \nexists$ \\ $0<\lambda <1 \Rightarrow$ ? \\ $\lambda = 1 \Rightarrow$ ?\\ $\lambda >1 \Rightarrow$ ?} & $\nexists$\\ \hline
         Product & $\nexists$ & $\nexists$ & $\nexists$ & $\nexists$ & $\nexists$ & $\nexists$ \\ \hline
         Schweizer-Sklar & $\nexists$ & $\nexists$ & \makecell{$\lambda<1 \Rightarrow \nexists$ \\ $\lambda=1 \Rightarrow IRD$ \\ $\lambda >1 \Rightarrow \exists CR$} & $\nexists$ & \makecell{$\lambda<1 \Rightarrow \nexists$ \\ $\lambda = 1 \Rightarrow$ ? \\ $\lambda > 1 \Rightarrow$ ? \\ $\lambda=+\infty \Rightarrow \nexists$} & $\nexists$\\ \hline
         Hamacher & $\nexists$ & $\nexists$ & $\nexists$ & $\nexists$ & $\nexists$ & $\nexists$ \\
         
         \hline \hline
         \begin{tabular}{@{}l@{}} Weak \\ decomposition \end{tabular} & $\nexists$ & $IDR$ & $\exists CR$ & $IDR$ & \makecell{$\lambda=-\infty \Rightarrow IDR$ \\ $-\infty < \lambda \le 0 \Rightarrow \nexists$ \\ $0<\lambda<+\infty \Rightarrow$ ?\\ $\lambda = +\infty \Rightarrow \nexists$}& \makecell{$\lambda <+\infty \Rightarrow IRD$ \\ $\lambda=+\infty \Rightarrow \nexists$} \\
         \hline
        \end{tabular}
        \caption{The representation of the relation between decomposition rules and decompositions. The symbol $\nexists$ means that there is no decomposition rule compatible with such decomposition, $\exists$CR means that there is a compatible decomposition rule with such decompositions, $IDR$ means that the decomposition induces a decomposition rule and ? means that the results in this article do not shed some light on the case. }
        \label{TablePreferences}
    \end{table}

\end{landscape}

\end{document}